\documentclass[12pt]{amsart}
\usepackage{amsmath,amssymb, xypic,verbatim,amscd,color}
\usepackage{graphicx}
\usepackage{colordvi}

\numberwithin{equation}{section}

\newcommand\frg{\mathfrak{g}}

\newcommand\mon{\overline{\operatorname{M}}_{0,n}}

\newtheorem{theorem}{Theorem}[section]
\newtheorem{remark}[theorem]{ Remark}

\newtheorem{proposition}[theorem]{Proposition}

\newtheorem{lemma}[theorem]{Lemma}

\begin{document}
\title[Level one conformal blocks divisors ]{Remarks on level one conformal blocks divisors}
%\author{Prakash Belkale}

\author{Swarnava Mukhopadhyay}
\address{Department of Mathematics\\ University of Maryland\\ CB \#4015, Mathematics Building
\\ College Park, MD 20742}
%\email{belkale@email.unc.edu}
\email{swarnava@umd.edu}
\begin{abstract} We show that conformal blocks divisors of type $B_r$ and $D_r$ at level one are effective sums of boundary divisors of $\mon$. We also prove that the conformal blocks divisor of type $B_r$ at level $1$ with weights $(\omega_1,\dots,\omega_1)$ scales linearly with the level. 
\end{abstract}
\maketitle
\section{Introduction}
Let $\frg$ be a simple Lie algebra, $\mathfrak{h}$ a Cartan subalgebra and $P_{\ell}(\frg)$, the set of dominant integral weights of $\frg$ at level $\ell$. Consider an $n$-tuple $\vec{\Lambda}=(\Lambda_1,\dots, \Lambda_n)$, where $\Lambda_i\in P_{\ell}(\frg)$. Corresponding to this data, consider conformal blocks bundles $\mathbb{V}_{\vec{\Lambda}}(\frg,\ell)$ on $\mon$. The first Chern classes of $\mathbb{V}_{\vec{\Lambda}}(\frg,\ell)$ are denoted by $\mathbb{D}(\vec{\Lambda},\frg,\ell)$. We refer the reader to  ~\cite{Sor} and ~\cite{TUY} for more details about conformal blocks. Due to the work of N. Fakhruddin, conformal blocks divisors play a central role in the birational geometry of $\mon$.

In this note, we focus on level one conformal blocks divisors of type $B_r$ and $D_r$. N. Fakhruddin has shown that conformal blocks divisors at level are often extremal in the Nef cone of $\mon$. Further it follows directly  from Chern class formulas of Fakhruddin (cf \cite{F}) that the conformal blocks divisors at level one for type $A_1$, $A_2$, $E_6$, $E_7$, $E_8$, $G_2$ and $F_4$ are effective sums of boundary divisors. We refer the reader to Sections 5.2.5-5.2.8 in \cite{F} for more details. It was recently proved in \cite{Mak} that all level one conformal divisors of type $A_r$ are effective combinations of boundary divisors. We prove the following:

\begin{theorem}\label{main}
The level one conformal blocks divisors of type $B_r$ and $D_r$ are effective combinations of boundary divisors. 
\end{theorem}
We do not know if conformal blocks of type $C_r$ at level one are effective combinations of boundary divisors. We believe that this is closely related to the behavior of conformal blocks of $A_1$ at level $r$ due to rank-level duality. It will be very interesting to find a non-boundary conformal blocks divisor on $\mon$. We now discuss the behavior of certain conformal blocks for type $B_r$ at level one under scaling. 

\begin{theorem}\label{secondary}
Let $n$ be an even integer and $\vec{\Lambda}$ be the $n$-tuple of weights $(\omega_1,\dots,\omega_1)$ of $B_r$ at level $1$. Then we have the following equality in $\operatorname{Pic}(\mon)$:
$$\mathbb{D}(N\vec{\Lambda},B_r,N)=N.\mathbb{D}(\vec{\Lambda}, B_r,1),$$ where $N\vec{\Lambda}$ is the $n$-tuple of level $n$ weights $(N\omega_1,\dots, N\omega_1)$. 
\end{theorem}

\subsection{Idea of Proof of Theorem \ref{main}}
Consider the complete graph $\Gamma[n]$ on $n$-vertices. Let $\Gamma[S]$ denote the set of edges of $\Gamma[n]$ and assume that vertices of $\Gamma[n]$ are labeled by the set $[n]=\{1,\dots, n\}$. To every edge $s$ of $\Gamma[n]$, we assign a rational number $w(s)$. We denote by $w(i)$, the total weight of all the edges passing through the vertex $i$. For a partition $I, J$ of $[n]$ with $|I|, |J|>2$, consider the set $V(I|J)$ of all vertices that starts at $I$ and ends in $J$. We denote the total weight of the edges in $V(I|J)$ by $w(I|J)$. 

Every divisor on $\mon$ can be written in the form $D=\sum_{i=1}^n a_i\psi_i-\sum_{I,J}c_{I,J} D_{I,J},$
where $\psi_i$ are the $i$-th psi classes and $D_{I,J}$ is a boundary divisor of $\mon$ corresponding to a partition $I, J$ of $[n]=\{1,\dots, n\}$. The following Lemma is proved in ~\cite{Mak}:
\begin{lemma}\label{maxlemma} A divisor $D$ is $\mathbb{Q}$-equivalent to an effective combinations of boundary divisors on $\mon$ if and only if there exists a $\mathbb{Q}$-valued weight function $w:\Gamma[S]\rightarrow \mathbb{Q}$ such that $w(i)=a_i$ and $w(I|J)$ is at least $c_{I,J}$ for all partitions $I, J$ of $[n]$. 
\end{lemma}
  In ~\cite{F}, Fakhruddin gave a formula for the first Chern class of a conformal blocks bundle on $\mon$ for arbitrary $\frg$ and level $\ell$. The author in ~\cite{M1} rewrote Fakhruddin's formula and expressed it as follows:
\begin{eqnarray*}
\mathbb{D}(\vec{\Lambda},\frg,\ell))&=& \operatorname{rk}\mathbb{V}_{\vec{\Lambda}}(\frg,\ell)  \bigg( \sum_{j=1}^n \Delta_{\Lambda_j}(\frg,\ell)\psi_j\bigg)\\
&& -\sum_{i=2}^{\lfloor \frac{n}{2} \rfloor}\epsilon_i  \sum_{\substack{I \subset \{1,2,\dots,n \} \\ |I|=i\\ J=[n]\setminus I\\}}\left\{\sum_{\Lambda \in P_{\ell}(\frg)}\Delta_{\Lambda}(\frg,\ell).\operatorname{rk}\mathbb{V}_{\vec{\Lambda}_I,\Lambda}(\frg,\ell).\operatorname{rk}\mathbb{V}_{{\vec{\Lambda}}_{J},\Lambda^*}(\frg,\ell)\right\}D_{I,J},
\end{eqnarray*}
 where $\vec{\Lambda}_I \subset \vec{\Lambda}$ denotes the set of weights $\Lambda_i$ such that $i \in I$, $\epsilon_i=\frac{1}{2}$ if $i=n/2$ and one otherwise. To complete the proof of Theorem \ref{main}, we construct a $\mathbb{Q}$-valued weight function $w$ on $\Gamma[S]$ satisfying the hypothesis of Lemma \ref{maxlemma}. In the next two sections, we give explicit constructions of weight functions $w$ for type $B_r$ and $D_r$ respectively. 
\begin{remark} 
%Theorem \ref{main} answers a question of Fakhruddin about limits of level one conformal blocks divisors of type $D_{r}$ being semi-ample. 
%We refer the reader to Remark 5.7 in ~\cite{F} for more details. 
N. Fakhruddin showed that conformal blocks bundles are globally generated and hence induces morphisms $\phi_{\mathbb{D}}$ from $\mon$.  It is interesting and challenging to classify the images of these morphisms. We hope that the weight functions constructed in the proof of Theorem \ref{main} may be used to identify the image of $\mon$ as ``Veronese quotients" and hopefully shed light on the cone of conformal blocks divisors in type $B_r$ and $D_r$. 
\end{remark}
\subsection{Acknowledgments} I thank P. Belkale, P. Brosnan and A. Gibney for discussions regarding this note.  After this work was done, I was informed by M. Fedorchuk about similar results he obtained about conformal blocks being effective combinations of boundary divisors.  %I also thank him for useful communications.
%\end{document}
\section{Proof of Theorem \ref{main} for $B_r$}The level one weights of $B_{r}$ are $\omega_0$, $\omega_1$ and $\omega_{r}$. We ignore $\omega_0$ completely due to ``Propagation of Vacua". The trace anomaly of the level one weights are $\Delta_{\omega_0}=0$, $\Delta_{\omega_1}=1/2$ and $\Delta_{\omega_{r}}=r(2r+1)/8.(2r)$. 

Let $n_1$, $n_2$ be the number of $\omega_1$'s, $\omega_{r}$'s in $\vec{\Lambda}$ respectively. If either $n_1$ or $n_2$ is zero, then the conformal blocks divisor is symmetric. Hence by a theorem of \cite{KM}, it is an effective combination of boundary divisors. If $n_2$ is odd, then the blocks is zero dimensional. Hence we assume that $n_1>0$ and $n_2=2m$ is a positive even number. The dimension of the corresponding conformal blocks at level one is $2^{m-1}$. 
\begin{remark}\label{factorb}
If $I, J$ be partition of $[n]$ and let the number of $\omega_{r}$'s in $I$ be even, then the conformal blocks factorizes into the direct sum $\mathbb{V}_{\vec{\Lambda}_I, \omega_0}(B_{r},1)\otimes \mathbb{V}_{\vec{\Lambda}_J,\omega_0}(B_{r},1) \oplus \mathbb{V}_{\vec{\Lambda}_I, \omega_1}(B_{r},1)\otimes \mathbb{V}_{\vec{\Lambda}_J,\omega_1}(B_{r},1)$ each of dimension $2^{m-2}$. On the other hand if the number of $\omega_{r}$'s are odd, it is isomorphic to $\mathbb{V}_{\vec{\Lambda}_I,\omega_{r}}(B_{r},1)\otimes \mathbb{V}_{\vec{\Lambda}_J,\omega_{r}}(B_{r},1)$.
 \end{remark}
\subsection{Weight function} We now describe the weight function $w(s)$ for type $B_r$ associated to the complete graph $\Gamma[n]$, whose vertices are marked by $\Lambda_i$. 
\begin{itemize}
\item  $w(s)=\frac{\Delta_{\omega_1}2^{m-1}}{(n_1-1)}-\frac{2^{m-1}}{n_1(n_1-1)}$, where $s$ is an edge joining two vertices labeled by $\omega_1$. 
\item $w(s)=\frac{\Delta_{\omega_r}2^{m-1}}{(n_2-1)}-\frac{2^{m-1}}{n_2(n_2-1)}$, where $s$ is an edge joining two vertices labeled by $\omega_{r}$.
\item $w(s)=\frac{2^{m-1}}{n_1n_2}$, otherwise. 
\end{itemize}
It is clear that the flow through every vertex is $2^{m-1}\Delta_{\Lambda_i}(B_{r},1)$. Let $I,J$ be a partition of the set $[n]=\{1,\dots, n\}$ and suppose $a_1,b_1$ be the number of $\omega_1$'s in $I$ and $J$ respectively. Further let $a_2=|I|-a_1$ and $b_2=|J|-b_1$. The total flow $w(I|J)$ through the partition $I, J$ is given by the following: 
\begin{equation}
\frac{a_1b_1}{n_1-1}({\Delta_{\omega_1}.2^{m-1}}-\frac{2^{m-1}}{n_1})+\frac{(a_1b_2+a_2b_1)2^{m-1}}{n_1n_2}+\frac{a_2b_2}{n_2-1}({\Delta_{\omega_{r}}.2^{m-1}}-\frac{2^{m-1}}{n_2})
\end{equation} The following proposition and Remark \ref{factorb} tell us that the function $w(s)$ satisfies the hypothesis of Lemma \ref{maxlemma}. 
\begin{proposition} The total flow $w(I|J)$ satisfies the following:
\begin{itemize}
\item $w(I|J) \geq \Delta_{\omega_{r}}.2^{m-1}$, when $a_2,b_2$ are odd. 
\item $w(I|J) \geq \Delta_{\omega_1}2^{m-2}$, when $a_2,b_2$ are even.
\end{itemize}

\end{proposition}
%$$w(I|J)=\Delta_{\omega_1}(B_r,1).\frac{a_1.b_1.2^{m-1}}{(n_1-1)}+ \Delta_{\omega_r}(B_r,1).\frac{a_2.b_2.2^{m-1}}{(n_2-1)}.$$ It follows directly from Remark \ref{factorb}, that the weight function $w(s)$ satisfies the hypothesis of Lemma \ref{maxlemma}. 
\section{Proof of Theorem \ref{main} for $D_r$}The level one weights of $D_{r}$ are $\omega_0$, $\omega_1$, $\omega_{r-1}$ and $\omega_{r}$. The trace anomaly of the weights are given as $\Delta_{\omega_1}=1/2$, $\Delta_{\omega_{r-1}}=\Delta_{\omega_{r}}=r/8$ and $\Delta_{\omega_{0}}=0$. We ignore $\omega_0$ due to ``Propagation of Vacua". 

First we observe that conformal blocks divisors of $D_{3}$ at level one are up to scaling same as conformal blocks divisor of $A_3$ at level one. These are all boundary by \cite{F} and \cite{Mak}. Hence assume that $r \geq 3$. Let $n_1$ be the number of $\omega_1$'s in $\vec{\Lambda}$ and $n_2=n-n_1$. As before we can assume that $n_1\neq 1$ and $n_2>1$. It follows from ~\cite{F}, that the level one conformal blocks of type $D_{r}$ with weights $\vec{\Lambda}$ is one dimensional if and only if $\sum_{i=1}^n\Lambda_i$ is in the root lattice of $D_r$ and zero otherwise.

\subsection{Weight function} We now describe the weight function $w(s)$ for type $B_r$ associated to the complete graph $\Gamma[n]$, whose vertices are marked by $\Lambda_i$. 
\begin{itemize}
\item  $w(s)=\frac{\Delta_{\omega_1}}{(n_1-1)}-\frac{1}{n_1(n_1-1)}$, where $s$ is an edge joining two vertices labeled by $\omega_1$. 
\item $w(s)=\frac{\Delta_{\omega_r}}{(n_2-1)}-\frac{1}{n_2(n_2-1)}$, where $s$ is an edge joining two vertices labeled either by $\omega_{r-1}$ or $\omega_{r}$.
\item $w(s)=\frac{1}{n_1n_2}$, otherwise. 
\end{itemize}It is clear that the flow through every vertex is $\Delta_{\Lambda_i}$. Let $I, J$ be a partition of the set $[n]=\{1,\dots, n\}$ and suppose $a_1,b_1$ be the number of $\omega_1$'s in $I$ and $J$ respectively. Further let $a_2=|I|-a_1$ and $b_2=|J|-b_1$. The total flow $w(I|J)$ through the partition $I,J$ is given by the following: 
\begin{equation}
\frac{a_1.b_1}{(n_1-1)}({\Delta_{\omega_1}}-\frac{1}{n_1})+ \frac{a_2.b_2}{(n_2-1)}({\Delta_{\omega_{r}}}-\frac{1}{n_2}) + \frac{a_1b_2+a_2b_1}{n_1n_2}
\end{equation}
Since the rank of the bundle $\mathbb{V}_{\vec{\Lambda}}(D_r,1)$ is one, it follows that for a partition $I, J$ of $[n]$, there is exactly one $\Lambda \in P_1(\frg)$ such that $\mathbb{V}_{\vec{\Lambda}_I,\Lambda}(D_r,1)$ and $\mathbb{V}_{\vec{\Lambda}_J,\Lambda^*}(D_r,1)$ are both of rank one. It is easy to see that $w(s)$ satisfies the hypothesis of Lemma \ref{maxlemma} from the following Proposition:
\begin{proposition} The total flow across a partition satisfies the following properties:
\begin{itemize}
\item $w(I|J) \geq \Delta_{\omega_1}$, if $a_2$ or $b_2$ is zero. 
\item $w(I|J) \geq \Delta_{\omega_{r}}$, otherwise.
\end{itemize}
\end{proposition}
%Let $\vec{\Lambda}$ denote an $n$-tuple of level one weights (of $D_{\ell}$ such the sum of the weights is in the root lattice of $D_{\ell}$. It follows from ~\cite{F}, that the level one conformal blocks of type $D_{\ell}$ with weights $\vec{\Lambda}$ is one dimensional. 

\section*{Proof of Theorem \ref{secondary}} Let $\sigma$ be the non trivial affine Dynkin diagram automorphism of type $B_r$. The automorphism $\sigma$ sends $N\omega_1$ to $N.\omega_0$, where $\omega_0$ is the zero-th affine fundamental weight of $B_r$. Since $n$ is an even integer, it follows from a result of ~\cite{FS} and ``Propagation of Vacua" that the rank of the conformal blocks bundle $\mathbb{V}_{N\vec{\Lambda}}(B_r,N)$ is one. One the other hand it is well known that the rank of the conformal blocks bundle $\mathbb{V}_{\vec{\Lambda}}(B_r,1)$ is also one if $n$ is even and zero otherwise. We refer the reader to \cite{F} for a proof. Now the proof of Theorem \ref{secondary} follows directly from induction on $N$ and applying Proposition 18.1 in \cite{BGM}. 
%Next we briefly recall the main idea behind the proof of Theorem \ref{main}.
\bibliographystyle{plain}

\begin{thebibliography}{10}
%\bibitem[AGSS]{agss}
%\bibitem{ABI}D. Altschuer, M. Bauer, C. Itzykson, {\em The branching rules of conformal embeddings}, Comm. Math. Phys. {\bf 132} (1990), no. 2, 349-364.
%\bibitem[AGS]{ags}

%\bibitem {A} T. Abe, {\em Strange duality for parabolic symplectic bundles on a pointed projective line}, Int. Math. Res. Not. IMRN { 2008}, Art. ID rnn121, 47 pp.

\bibitem {BGM} P. Belkale, A. Gibney, S. Mukhopadhyay, {\em Quantum cohomology and conformal blocks on $\mon$}, arXiv:1308:4906v2.
%\bibitem {BB} A. Bais, P. Bouwknegt, {\em A classification of subgroup truncations of the bosonic string}, Nuclear Phys. B {\bf 279}(1987), no. 3-4, 561-70.
%\bibitem {BF} E. Frenkel, D. Ben-Zvi, {\em Vertex Algebras on algebraic curves}, Mathematical Surveys and Monographs {\bf 88}, AMS 2001.
%\bibitem {BP}  A. Boysal, C. Pauly, {\em Strange duality for Verlinde spaces of exceptional groups
%at level one}, Int. Math. Res. Not. 2009; doi:10.1093/imrn/rnp151.
%\bibitem{FG} G. Farkas, A. Gibney, {\em Mori cones of moduli spaces of pointed curves of small genus }, Trans. Amer. Math. Soc. {\bf 355} (2003), no. 3, 1183-1199 (electronic), DOI 10.1090/S0002-9947-02-03165-3.
\bibitem{FS}  J. Fuchs, C. Schweigert, {\em The Action of Outer Automorphisms on Bundles of Chiral Blocks}, Comm. Math. Phys. 206 (1999), no. 3, 691-736.
\bibitem{F}  N. Fakhruddin,{\em Chern classes of conformal blocks,} Contemp. Math., {\bf 564}, Amer. Math. Soc., Providence, RI, 2012, 145-176.
\bibitem{KM} S. Keel, J. McKernan,{\em Contractible extremal rays on $\mon$}, arXiv:alg-geom/9607009v1.
\bibitem{Mak} M. Fedorchuk, {\em New Nef divisors on $\mon$}, arXiv:1308.5593v1.
%\bibitem[H]{H} K. Hasegawa,
%\bibitem[GG]{GG}
%\bibitem {KW} V. Kac, M. Wakimoto, {\em Modular and conformal invariant constraints in representation theory of affine algebras}, Advances in Mathematics {\bf 70} (1988):156-234.
\bibitem {M1} S. Mukhopadhyay, {\em Rank-level duality and conformal block divisors}, arXiv:1308.0854v1.
%\bibitem {M2}S. Mukhopadhyay,  {\em Diagram automorphisms to rank-level duality}, arXiv:1308:1756v1.
%\bibitem {NT} T. Nakanishi, A. Tsuchiya, {\em Level-rank duality of WZW models in conformal field theory}, Comm. Math. Phys. Volume {\bf 144} (1992), no. 2, 351-372.
%\bibitem[2]{ABI}D. Altschuer, M. Bauer, C. Itzykson, {\em The branching rules of conformal embeddings}, Comm. Math. Phys. 132 (1990), no. 2, 349-364. 
%\bibitem[9]{TK} A. Tsuchiya, Y.  Kanie, {\em  Vertex operators in conformal field theory on $\mathbb{P}^1$ and monodromy representations
%of braid group}, Conformal field theory and solvable lattice models (Kyoto, 1986), 297-372, Adv. Stud. Pure Math., 16, Academic Press, Boston, MA, 1988.
%\bibitem[9]{KW}
\bibitem {Sor} C. Sorger, {\em Le formula de Verlinde}, Ast{$\acute{e}$}rique, (1996), pp. Exp. No. 794, 3, 87-114. S$\acute{e}$minaire Bourbaki, Vol. 1994/95.
\bibitem{TUY} A. Tsuchiya, K. Ueno, Y. Yamada, {\em Conformal field theory on universal family of stable curves with
gauge symmetries}, Integrable systems in quantum field theory and statistical mechanics, 459-566,
%Adv. Stud. Pure Math. {\bf  19}, Academic Press, Boston, MA, 1989.
%\bibitem {T} Y. Tsuchimoto, {On the coordinate-free description of the conformal blocks}, J. Math. Kyoto Univ. {\bf 33} (1993) 29-49.
%\bibitem[26]{Ueno} K. Ueno, {\em Conformal field theory with gauge symmetry}, Fields Institute Monographs, 24. American Mathematical Society, Pr


\end{thebibliography}
\def\noopsort#1{}
%\begin{thebibliography}{10}

%\end{document}
\end{document}